\documentclass[11pt]{article}

\usepackage{amsmath}
\usepackage{amssymb}
\usepackage{eucal}

\begin{document}

\baselineskip 16pt

\title{On one generalization of modular subgroups \thanks{Research
 is supported by an NNSF grant of China (Grant No.
11401264) and a TAPP of Jiangsu Higher Education Institutions (PPZY 2015A013)}}

\author{Jianhong Huang\\\
{\small School of Mathematics and Statistics, Jiangsu Normal University,}\\
 {\small Xuzhou, 221116, P.R. China}\\
{\small E-mail: jhh320@126.com}\\ \\
Bin Hu \thanks{Corresponding author} \\
{\small School of Mathematics and Statistics, Jiangsu Normal University,}\\
{\small Xuzhou 221116, P. R. China}\\
{\small E-mail: hubin118@126.com}\\ \\
{Xun Zheng}\\
{\small School of Mathematics and Statistics, Jiangsu Normal University,}\\
{\small Xuzhou 221116, P. R. China}\\
{\small E-mail: zx309556477@126.com }}

\date{}
\maketitle

\begin{abstract}  Let $G$ be a finite group.
 If  $M_n < M_{n-1} < \ldots < M_1 < M_{0}=G  $ where $M_i$ is a
 maximal subgroup of  $M_{i-1}$ for all $i=1, \ldots ,n$, then $M_n $ ($n > 0$)
 is an \emph{$n$-maximal subgroup}  of $G$.
 A subgroup $M$ of
 $G$ is called  \emph{modular}  if the following conditions are held:
(i) $\langle X, M \cap Z \rangle=\langle X, M \rangle \cap Z$ for all $X \leq G, Z \leq
 G$ such that
$X \leq Z$, and
(ii) $\langle M, Y \cap Z \rangle=\langle M, Y \rangle \cap Z$ for all $Y \leq G, Z \leq
 G$ such that
$M \leq Z$.

In this paper, we  study finite groups whose   $n$-maximal subgroups are
   modular.

\end{abstract}

\footnotetext{Keywords: finite group,  modular subgroup, $n$-maximal subgroup,
 nearly nilpotent group, strongly supersoluble group.}

\footnotetext{Mathematics Subject Classification (2010): 20D10,
20D15, 20D20}
\let\thefootnote\thefootnoteorig

\section{Introduction}

Throughout this paper, all groups are finite and $G$ always denotes
a finite group.  Moreover,  $\mathbb{P}$ is the set of all  primes and the
 symbol $\pi(G)$ stands for the set of prime divisors of
the order of $G$.

We say that $G$  is:  \emph{nearly nilpotent} if $G$ is supersoluble and $G$ induces
on any  its  non-Frattini chief factor $H/K$ (that is, $H/K\nleq \Phi (G/K)$)
 an  automorphism group of
order dividing a prime;
  \emph{srtongly supersoluble} if $G$ is
 supersoluble and $G$ induces
on any  its chief factor $H/K$ an  automorphism group of square free
order.  We use    $\frak{N}_{n}$  and  $\frak{U}_{s}$ to denote the
classes  of all nearly nilpotent and
of all strongly  supersoluble groups, respectively. Nearly nilpotent and strongly supersoluble groups were
studied  respectively in \cite{1-3}  and  \cite{Zim, Vas}.

 It is clear that:  the group $C_{7}\rtimes \text{Aut}(C_{7})$ is strongly
supersoluble but it is not nearly nilpotent; the group $C_{13}\rtimes \text{Aut}(C_{13})$ is
supersoluble but  it is not strongly supersoluble; the group $S_{3}$ is
 nearly nilpotent but  it is not  nilpotent.

A  subgroup $M$ of  $G$ is called \emph{ modular} if $M$ is a modular element
(in the sense of
 Kurosh \cite[ p. 43]{Schm})  of the
 lattice ${\cal L}(G)$ of all subgroups of $G$, that is,

(i) $\langle X,M \cap Z \rangle=\langle X, M \rangle \cap Z$ for all $X \leq G, Z \leq
 G$ such that
$X \leq Z$, and

(ii) $\langle M, Y \cap Z \rangle=\langle M, Y \rangle \cap Z$ for all $Y \leq G, Z \leq
 G$ such that
$M \leq Z$.

Recall
 that a subgroup $H$ of $G$ is called a \emph{ 2-maximal} (\emph{second maximal}) subgroup of G whenever $H$ is
a maximal subgroup of some maximal subgroup $M$ of $G$. Similarly we can define
 \emph{ 3-maximal} subgroups, and so on.

The relationship between $n$-maximal subgroups (where  $n>1$) of  $G$ and the
 structure of $G$ was studied by many authors (see, in particular, the recent  papers
\cite{li1}--\cite{K2} and Chapter 4 in the
book \cite{GuoII}).  One of the earliest results
in this   line research was obtained  by  Huppert in the article \cite{HupI} who established the supersolubility
 of the group whose
all second maximal subgroups are normal. In the same article
 Huppert proved that if all 3-maximal
subgroups of $G$ are normal in $G$, then the commutator subgroup $G'$
 of $G$ is a nilpotent group and the
principal rank of $G$ is at most 2.
  These  two  results  were  developed  by  many  authors.
In particular,  Schmidt proved  \cite{1-3} that: if all 2-maximal
subgroups of $G$ are modular in $G$, then $G$ is nearly nilpotent;  if all 3-maximal subgroups
 of $G$  are  modular in $G$ and $G$ is not supersoluble, then  either
    $G$ is
a group of order $pq^{2}$ for primes $p$ and $q$ or $G=Q\rtimes P$,
 where $Q=C_{G}(Q)$ is a quaternion group of order 8 and $|P|=3$.
Mann proved \cite{mann} that if all $n$-maximal subgroups
 of a soluble group  $G$
are subnormal and
$n < |\pi (G)|$,
then $G$ is nilpotent; but if  $n \leq |\pi (G)| + 1$,
then $G$ is  $\phi$-dispersive for some ordering  $\phi$ of
$\mathbb P$.  Finally, in the case  $n\leq |\pi (G)|$
Mann described $G$ completely.

 In this paper, we  prove the following modular analogues of the above-mentioned
 Mann's results.

 {\bf Theorem A.}  {\sl  Suppose that $G$ is soluble and every $ n$-maximal subgroup
 of $ G$  is  modular.  If $n \leq |\pi (G)|$, then
  $G$ is  strongly supersoluble and $G$ induces
on any  its  non-Frattini chief factor $H/K$  an  automorphism group of order $p_{1}\cdots p_{m}$ where $m\leq n$
 and $p_{1},  \ldots
 , p_{m}$ are distinct primes.   }

 We use $G^{\frak{U}_{s}}$ to denote the intersection of all normal
 subgroups $N$ of $G$
with  strongly supersoluble quotient $G/N$.

 {\bf Theorem B.}  {\sl  Suppose that $G$ is soluble and every $ n$-maximal subgroup
 of $ G$
 is modular.  If $n \leq |\pi (G)|+1$, then $G^{\frak{U}_{s}}$ is a nilpotent Hall subgroup of $G$.  }

Finally, note that the restrictions on $|\pi (G)|$  in Theorems A  and
 B cannot be weakened (see Section 4 below).

\section{Proof of theorem A}

A normal subgroup $A$ of  $G$ is said to be \emph{hypercyclically embedded} in
 $G$  \cite[p. 217]{Schm} if either $A=1$ or $A\ne 1$ and
every chief factor of $G$ below $A$ is cyclic.
We use $Z_{\mathfrak{U}}(G)$ to denote the product of all  normal
hypercyclically embedded subgroups of $G$. It is clear that a normal
subgroup $A$ of $G$ is hypercyclically embedded in $G$ if and only
if $A\leq Z_{\mathfrak{U}}(G)$.

Recall that $G$ is said to be a \emph{$P$-group} \cite[p. 49]{Schm}
if $G=A\rtimes \langle t \rangle$ with an elementary abelian
$p$-group $A$ and an element $t$ of
 prime order $q\ne p$ induces a non-trivial power
 automorphism  on $A$.

The following two lemmas collect the properties of modular subgroups which
we use in our proofs.

{\bf Lemma 2.1}  (See Theorems 5.1.14 and 5.2.5 in \cite{Schm}).
 {\sl Let $M$  be a modular subgroup of $G$.}

(i) {\sl  $M/M_{G}$ is nilpotent  and $M^G/M_G \leq Z_{\mathfrak{U}}(G/M_G)$}.

(ii) {\sl If    $M_{G}=1$, then
 $$ G=S_{1}\times \cdots \times S_{r} \times K,$$
 where  $0\leq r\in \mathbb{Z}$ and for all $i, j \in \{1, \ldots , r\}$,
}

(a) {\sl $S_{i}$ is a non-abelian $P$-group,}

(b) {\sl  $(|S_{i}|, |S_{j}|)=1=(|S_{i}|, |K|)$ for all $i\ne j$, }

(c) {\sl $M=Q_{1}\times \cdots \times Q_{r}\times (M\cap K)$ and $Q_{i}$
 is a non-normal Sylow subgroup of $S_{i}$, }

(d) {\sl $M\cap K$ is quasinormal in $G$. }

{\bf Lemma 2.2} (See  p. 201 in \cite{Schm}). {\sl Let $A$, $B$  and
$N$ be subgroups of $G$, where $A$ is modular in $G$ and $N$ is
normal in $G$.}

(1) {\sl  If  $B$ is modular in $G$, then  $\langle
A,  B\rangle $  is  modular in $G$.}

(2) {\sl  $AN/N$
 is modular in  $G/N$}.

(3) {\sl   $N$ is   modular in $G$.}

(4) {\sl If $A \leq B$, then  $A$  is modular in  $B$.}

(5) {\sl If $\varphi$ is an isomorphism of $G$ onto $\bar{G}$, then
 $A^{\varphi}$ is modular in $\bar{G}$.}

A subgroup $H$ of $G$ is said to be  \emph{quasinormal }  (respectively
 \emph{$S$-quasinormal }) in $G$ if
$HP=PH$ for all subgroups  (for all Sylow subgroups) $P$ of $G$.

{\bf Lemma 2.3} (See Chapter 1 in  \cite{prod}).  {\sl Let $H\leq K\leq
G$.  }

(1) {\sl If $H$ is $S$-quasinormal   in $G$, then $H$ is $S$-quasinormal
in $K$.  }

(2) {\sl  Suppose that $H$ is normal in $G$. Then $K/H$ is
$S$-quasinormal   in $G/H$ if and only if $K$ is $S$-quasinormal
in $G$.   }

(3) {\sl If $H$ is $S$-quasinormal   in $G$, then $H$ is subnormal in $G$ and
 $H^{G}/H_{G}$ is nilpotent.}

{\bf Lemma 2.4.}  {\sl Suppose that   $G$ is  soluble,  and let $N\ne G$ be a minimal
normal subgroup of $G$. Suppose also that every $ n$-maximal subgroup
 of $ G$
 is either modular or $S$-quasinormal in $G$, where  $n\leq|\pi (G)|+ r $ for some integer $r$.
 Then there is
 a natural number $m\leq n$ such that  every $ m$-maximal subgroup
 of $ G/N$  is either modular or $S$-quasinormal in $G/N$ and $m \leq |\pi (G/N)|+ r $.}

{\bf Proof. }  First assume that $N$ is not a Sylow subgroup of $G$. Then
 $|\pi (G/N)|=|\pi (G)|$.
Moreover, if $H/N$ is an $n$-maximal subgroup of $G/N$,  then $H$ is
an $ n$-maximal subgroup of $G$, so  $H$ is either modular or
$S$-quasinormal in $G$ by hypothesis.
 Consequently, $H/N$ is either modular or $S$-quasinormal in $G/N$
by Lemmas 2.2(2) and 2.3(2). On the other hand, if $G/N $includes no
$n$-maximal subgroups,  then, by the solubility
 of $G$, the trivial
subgroup of $G/N$ is modular in $G/N$ and is a unique $m$-maximal
subgroup of $G/N$ for some
 $m < n$ with $m < |\pi(G/N)|$.   Hence  $m < |\pi (G/N)|+ r $.
 Thus the conclusion of the lemma is fulfilled for $G/N$.

 Finally,
consider the case that $N$ is a Sylow $p$-subgroup of $G$.
 Let $E$ be a Hall $p'$-subgroup of $G$.
 It is clear that $|\pi (E)| = |\pi (G)| -1$   and  $E$ is a maximal subgroup   of $G$.
 Therefore,
every  $(n - 1)$-maximal subgroup of $E$ is either modular or $S$-quasinormal
 in $E$ by Lemmas 2.2(4) and 2.3(1).
 Thus, by the isomorphism $G/N\simeq E$, Lemma 2.2(5) implies that
 every  $(n - 1)$-maximal subgroup
of $G/N$ is either modular or $S$-quasinormal  in $G/N$,
 and also we have $n-1 \leq |\pi (G/N)|+ r $. The lemma is proved.

A \emph{formation} is a class $\frak{F}$ of groups with the
following properties: (i) Every homomorphic image of any group in
$\frak{F}$ belongs to $\frak{F}$; (ii) $G/N\cap R\in \frak{F}$ whenever $G/N \in \frak{F}$ and
 $G/R \in \frak{F}$. A formation $\frak{F}$ is said to be:
\emph{saturated} if $G\in \frak{F}$ whenever $G/\Phi (G) \in \frak{F}$;
 \emph{hereditary} if $H\in \frak{F}$ whenever $H\leq G \in
\frak{F}$.

{\bf Lemma 2.5}  (See Theorem A in \cite{Vas}).  {\sl  The class
 of all strongly supersoluble groups is a hereditary saturated formation. }

Let $\frak{X}$ be a  class  of groups.  A group $G$ is called a
\emph{minimal non-$\frak{X}$-group} \cite{GuoII} or
\emph{$\frak{X}$-critical group}  \cite{DH} if $G$ is not in
$\frak{X}$ but all proper subgroups of $G$ are in $\frak{X}$.  An
$\frak{N}$-critical group is  also
 called  a \emph{Schmidt group}.

Fix some ordering $\phi$ of $\mathbb{P}$. The record $p\phi q$ means that $p$ precedes $q$
 in $\phi$ and
 $p \ne q$. Recall that  a group $G$ of order $p_1^{\alpha _1}p_2^{\alpha _2}\dots
p_n^{\alpha _n}$ is called \emph{$\phi$-dispersive} whenever
$p_1\phi p_2\phi \dots \phi p_n$  and  for every $i$ there is a
normal subgroup of $G$ of order $p_1^{\alpha _1}p_2^{\alpha _2}\dots
p_i^{\alpha _i}$.  Furthermore, if $\phi$  is such that $p\phi  q$
always implies $p > q$ then every $\phi$-dispersive group is called
\emph{Ore dispersive}.

{\bf Lemma 2.6} (See \cite[I, Propositions 1.8, 1.11 and 1.12]{GuoII}). {\sl The following
 claims hold for every  $\frak{U}$-critical group  $G$:}

(1) {\sl $G$ is soluble and $|\pi(G)| \leq  3$.}

(2) {\sl If G is not a Schmidt group, then G is Ore dispersive.}

(3) {\sl $ G^{\frak{U}}$ is a unique normal Sylow subgroup of $G$.}

(4) {\sl    If $S$ is a complement to $ G^{\frak{U}}$ in $G$, then $S/S \cap \Phi (G)$ is
 either a cyclic prime power order group or a  Miller-Moreno
 (that is, a minimal non-abelian) group.  }

(5) {\sl $G^{\mathfrak{U}}/\Phi (G^{\mathfrak{U}})$ is a non-cyclic chief factor of $G$.}

(6) {\sl If $G^{\mathfrak{U}}$ is non-abelian,  then the center, commutator subgroup,
 and Frattini subgroup of $G^{\frak{U}}$ coincide with one another. }

(7) {\sl If  $p>2$, then $G^{\mathfrak{U}}$ is of exponent  $p$; for $p =
2$
 the exponent of $G^{\mathfrak{U}}$ is at most 4. }

{\bf Lemma 2.7} (See Lemma 12.8 in \cite{1989}). {\sl If $H/K$ is an
abelian chief factor of $G$ and $M$ is a maximal subgroup of $G$
such that $K \leq  M$ and $MH = G$, then
$$G/M_{G} \simeq (H/K)\rtimes (G/C_{G}(H/K))\simeq (HM_{G}/M_{G})\rtimes
 (G/C_{G}(HM_{G}/M_{G})) .$$}

 The following lemma is evident.

{\bf Lemma 2.8.}   {\sl
 If $H/K$ and $T/L$ are $G$-isomorphic chief factors of $G$, then
 $(H/K)\rtimes (G/C_{G}(H/K))\simeq  (T/L)\rtimes (G/C_{G}(T/L))$}.

Recall that a class of soluble groups $\frak{X}$ is a
 \emph{Schunck class} \cite[III, 2.7]{DH} if
$G\in \frak{X}$ whenever $G/M_{G}\in\frak{X}$ for all maximal subgroups
$M$ of $G$.

{\bf Proposition 2.9.} {\sl  The
 class of all nearly nilpotent  groups ${\frak{N}}_{n}$ is a  Schunck class, and
 ${\frak{N}}_{n} \subseteq {\frak{U}}_{s}$. Hence every homomorphic image of
 any nearly nilpotent group  is   nearly nilpotent, and $G$ is  nearly
nilpotent whenever $G/\Phi (G)$ is  nearly nilpotent.}

{\bf Proof.}  Suppose that for every maximal subgroup $M$ of $G$ we have
$G/M_{G} \in  {\frak{N}}_{n}$.  Then $G/\Phi (G)$ is supersoluble, so $G$
is supersoluble.
 If $H/K$ is a non-Frattini
chief factor of $G$ and $M$ is a maximal subgroup of $G$ such that
$K \leq  M$ and $MH = G$, then $G/M_{G} \simeq (H/K)\rtimes
(G/C_{G}(H/K))$ by Lemma 2.7. Since clearly
 $C_{(H/K)\rtimes (G/C_{G}(H/K))}(H/K)=H/K$, it follows that
$|G/C_{G}(H/K)|=p$ is a prime. Hence  $G\in {\frak{N}}_{n}$. Therefore
${\frak{N}}_{n}$ is a  Schunck class, so  every homomorphic image of
 any nearly nilpotent group  is   nearly nilpotent, and $G$ is  nearly
nilpotent whenever $G/\Phi (G)$ is  nearly nilpotent by  \cite[III, 2.7]{DH} .

Now  we show that every   nearly nilpotent  group $G$ is strongly
supersoluble. Assume that this is false and let $G$ be a
counterexample of minimal order. Let $R$ be a minimal normal
subgroup of $G$.  Then $G/R$ is strongly supersoluble by the choice
of $G$ since $G/R$ is nearly nilpotent. Moreover, if $R\leq \Phi
(G)$, then  $G$ is strongly supersoluble by Lemma 2.5, contrary to
the choice of $G$.
 Therefore   $R\nleq \Phi (G)$, so $G/C_{G}(R)$ is of  prime order since $G$ is
nearly nilpotent. Therefore $G$ is strongly supersoluble by the Jordan-H\"{o}lder theorem.
 This
contradiction completes the proof of the proposition.

{\bf Lemma 2.10.}  {\sl Let $G=R\rtimes M$ be a soluble primitive group, where
$R=C_{G}(R)$ is a minimal normal subgroup of $G$. Let $T\ne 1$ be a subgroup of $G$.
 Suppose that $G$ is not nearly nilpotent.}

 (1) {\sl If $T < M$, then $T$ is neither modular nor $S$-quasinormal in $G$.  }

 (2) {\sl If $T < R$ and $|M|$ is a prime, then  some subgroup $V$ of $R$ with $|V|=|T|$
 is neither modular nor $S$-quasinormal in $G$. }

{\bf Proof. }    (1)  First assume that
 $T$ is modular   in $G$ but it is not $S$-quasinormal in $G$.
 Then $T$ is not quasinormal in $G$, so  Lemma 2.1(ii)
 implies that $G$  is a non-abelian $P$-group since  $T_{G}\leq M_{G}=1$.
But then $G$ is  supersoluble.
This contradiction shows that $T$ is $S$-quasinormal in $G$, so $T$ is
subnormal in $G$ by Lemma 2.3(3). Hence $1 < T^{G}=T^{RM}=T^{M}\leq
M_{G}=1$ by \cite[A, 14.3]{DH}, a contradiction.    Hence we have (1).

(2)  Let $V$ be a subgroup of $R$ with $|V|=|T|$ such that $V$ is
normal in a Sylow $p$-subgroup of $G$.
 If $V$ is $S$-quasinormal in $G$, then for every Sylow $q$-subgroup $Q$ of $G$,
 where $q\ne p$, we have $VQ=QV$ and so $V=R\cap VQ$.
Hence $Q\leq N_{G}(V)$. Thus $V$ is normal in $G$, a contradiction. Hence
$V$ is modular in $G$, which implies that $1  < V\leq R\cap
Z_{\frak{U}}(G)$ by Lemma 2.1(i)  and so $R\leq Z_{\frak{U}}(G)$.
 But then $|R|=p$, which implies that  $G$ is nearly nilpotent, a contradiction.

The lemma is proved.

{\bf Proposition 2.11.}   {\sl  If    every  maximal subgroup of $G$
or every
  2-maximal subgroup of $G$  is either modular
  or $S$-quasinormal in $G$, then $G$ is nearly nilpotent. Hence
 $G$ is strongly supersoluble.}

{\bf Proof.}
 Assume this proposition is false and let $G$ be
 a  counterexample of minimal
order.

First we show that $G$ is soluble. Indeed, if $M$ is a maximal
subgroup of $G$ and either $M$ is modular in $G$ or $M$ is
$S$-quasinormal in $G$, then   $|G:M|$ is a prime  by Lemmas 2.1(i)
and 2.3(3). Therefore if every  maximal subgroup  of $G$  is either
modular
  or $S$-quasinormal in $G$, then $G$ is supersoluble. On the other hand,
if every
  2-maximal subgroup of $G$  is either modular
  or $S$-quasinormal in $G$, then every maximal subgroup of $G$ is
supersoluble by Lemmas 2.2(4) and 2.3(1) and so $G$ is soluble by Lemma 2.6(1).

Therefore, in view of Proposition 2.9, we need only to show that for
every maximal subgroup $M$ of $G$ we have $G/M_{G}\in {\frak{N}}_{n}$.
If $M_{G}\ne 1$, then the choice of $G$ and Lemmas 2.2(2) and 2.3(2) imply that
 $G/M_{G}\in {\frak{N}}_{n}$. Now assume that  $M_{G}=1$, so   there is a
minimal normal subgroup $R$ of $G$ such that $G=R\rtimes M$ and
$R=C_{G}(R)$ by \cite[A, 15.6]{DH}. Then  $M$ is  not $S$-quasinormal in
$G$ by Lemma 2.3(3). On the other hand, if $M$ is modular in $G$, then $G=M^{G}$
 is a non-abelian
$P$-group by Lemma 2.1(ii). It
follows that $G$ is nearly nilpotent, a contradiction.  Hence every
2-maximal subgroup of $G$ is either modular
  or $S$-quasinormal in $G$.

Now let $T$ be any  maximal subgroup of $M$. Then $T$  is either modular
  or $S$-quasinormal in $G$, so  $T=1$ and  hence  $|M|=q$ for some prime $q$. Therefore  $R$ is a
maximal subgroup of $G$.  Then every maximal subgroup of $R$ is
either modular or $S$-quasinormal in $G$ and  so  $|R|=p$ by Lemma 2.10(2),
 which implies that $|G|=pq$.
Hence $G$ is nearly nilpotent, a contradiction.

The proposition is proved.

In fact, Theorem A is a special case of the following

 {\bf Theorem 2.12.}  {\sl  Suppose that $G$ is soluble and every $ n$-maximal subgroup
 of $ G$
 is either  modular  or $S$-quasinormal in $G$.  If $n \leq |\pi (G)|$, then
  $G$ is  strongly supersoluble  and $G$ induces
on any  its  non-Frattini chief factor $H/K$
 an  automorphism group of order $p_{1}\cdots p_{m}$, where $m\leq n$
 and $p_{1},  \ldots
 , p_{m}$ are distinct primes.   }

{\bf Proof.}  Assume this theorem is false and let $G$ be
 a  counterexample  of minimal
order.

First we show that $G$ is strongly supersoluble. Suppose that this is
false.
Let $R$ be a minimal normal subgroup of $G$.

(1) {\sl $G/R$ is strongly supersoluble. Hence $G$ is primitive and
so $R\nleq \Phi (G)$ and $R=C_{G}(R)=O_{p}(G)$
 for some prime $p$. }

Lemma 2.4 implies that the hypothesis holds for $G/R$, so the choice
of $G$ implies that $G/R$ is strongly supersoluble. Therefore, again
by the choice of $G$, $R$ is a unique minimal normal subgroup of $G$
and  $R\nleq \Phi (G)$ by Lemma 2.5. Hence $G$ is primitive and so
$R=C_{G}(R)=O_{p}(G)$
 for some prime $p$ by \cite[A, 15.6]{DH}.

(2) {\sl Every maximal subgroup $M$ of $G$ is strongly supersoluble. }

By hypothesis  every $(n-1)$-maximal subgroup $T$ of $M$ is either
 modular or  $S$-quasinormal in $G$. Hence  $T$ is   modular
  in $M$ by Lemma 2.2(4)  in
 the former
case, and it is   $S$-quasinormal in $M$ by Lemma 2.3(1) in the second case.
 Since the solubility of $G$ implies          that either
  $|\pi (M)|=|\pi (G)|$  or $|\pi (M)| = |\pi (G)| - 1$, the
hypothesis holds for $M$.
 It follows that $M$ is  strongly supersoluble by the choice
of $G$.

(3) {\sl $G$ is supersoluble. }

Suppose that this is false. Since every maximal subgroup $M$ of $G$ is
 strongly supersoluble  by Claim
(2),   $G$ is a minimal non-supersoluble group. Then Lemma 2.6(1) yields that
  $ |\pi (G)| = 2$ or $ |\pi (G)| = 3$. But in the former case $G$ is
strongly
supersoluble by Proposition 2.11, so  $ |\pi (G)| = 3$ and every 3-maximal
subgroup of $G$ is either modular or $S$-quasinormal in $G$.
  Claim (1)  and
Lemma 2.6   imply  that $G=R\rtimes S$, where $S$ is a
Miller-Moreno group. Moreover, since $|\pi (S)|=2$ and $S$ is
strongly supersoluble,
 $S$ is not nilpotent and so $S =Q\rtimes T$,
 where
$|Q|=  q$, $|T|=t$ and    $C_{S}(Q)=Q$ for some distinct primes $q$ and
$t$ by \cite[I, Proposition 1.9]{GuoII}. Hence  $R$ is
 a 2-maximal subgroup of $G$, so every maximal subgroup of $R$ is either
 modular or $S$-qusinormal $G$. Therefore $G$ is supersoluble by Lemma
2.10(2).

(4) {\sl $G$ is strongly supersoluble.}

  From Claims (1) and (3) we get that for some maximal subgroup $M$ of $G$
we have
$G=R\rtimes M=C_{G}(R)\rtimes M$ and $|R|=p$, so $M$ is cyclic.
Since $G$  is not strongly supersoluble,  for some prime $q$ dividing
$|M|$ and for the Sylow $q$-subgroup $Q$ of $M$ we have $|Q| > q$. First
assume that $RQ\ne G$, and let $RQ\leq V$, where  $V$ is a maximal
subgroup of $G$. Then  $V$  is strongly supersoluble by Claim (2). Hence
$C_{Q}(R)\ne 1$, contrary to $R=C_{G}(R)$.  Hence  $RQ= G$ and   so $|\pi
(G)|=2$. Therefore $G$ is strongly supersoluble by Proposition 2.11, a
contradiction. Thus we have  (4).

(5) {\sl   $G$ induces
on any  its  non-Frattini chief factor $H/K$
 an  automorphism group $G/C_{G}(H/K)$ of order $p_{1}\cdots p_{m}$ where $m \leq n$
 and $p_{1},  \ldots
 , p_{m}$ are distinct primes.  }

If $G$ is nearly nilpotent, it is clear. Now suppose that $G$ is not nearly nilpotent.  Let $M$ be a maximal subgroup of $G$ such that $K\leq M$ and $MH=G$. Then
$G/M_{G} \simeq (H/K)\rtimes (G/C_{G}(H/K))$ by Lemma 2.7.
If $M_{G}\ne 1$, the choice of $G$ implies that $m\leq n$.
 Now suppose that $M_{G}=1$, so $G=H\rtimes
M$, where $|H|$ is a prime and $H=C_{G}(H)$. Then, by Claim (4), $M$ is a cyclic
group of order $p_{1} \ldots  p_{m}$ for some distinct primes  $p_{1},
\ldots , p_{m}$. Assume that $n < m$.
 Then $G$ has an $n$-maximal subgroup $T$ such that
$T\leq M$ and $|T|$ is not a prime. But since $G$ is not nearly nilpotent, this  is not
 possible by Lemma 2.10(1). This contradiction completes the proof of
the result.

\section{Proof  of Theorem  B}

{\bf Lemma 3.1}  (See  p. 359 in  \cite{DH}).  {\sl Given any
 ordering $\phi$ of the set of all primes, the class of all $\phi$-dispersive groups
 is a saturated  formation.}

{\bf Proposition 3.2.}   {\sl Suppose that
 every 3-maximal subgroup of $G$ is either $S$-quasinormal or modular in
$G$.   If $G$ is not  supersoluble, then either $G$ is
a group of order $pq^{2}$ for some distinct primes $p$ and $q$, or $G=Q\rtimes P$,
 where $Q=C_{G}(Q)$ is a quaternion group of order 8 and $|P|=3$.}

{\bf Proof.}  Assume that this proposition is   false and let $G$ be
 a  counterexample of minimal    order.  Lemmas 2.2(4), 2.3(1) and Proposition
2.11 imply that  every maximal subgroup of $G$ is strongly
supersoluble.  Hence $G$ is soluble by Lemma 2.6(1), so $|\pi (G)|=2$ by
Theorem  2.12.

Since   $G$ is not supersoluble, $G$  is a  $\frak{U}$-critical group.
 Let $D=G^{\frak{U}}$  be the
supersoluble residual of $G$. Lemma 2.6
 implies that   the following hold:
(a) $D$ is a Sylow $p$-subgroup of $G$ for some prime $p$, and if
 $Q$ is a Sylow $q$-subgroup of $G$, where $q\ne p$, then
 $DQ=G$ and $Q/Q \cap \Phi (G)$ is
 either a cyclic prime power order group or a  Miller-Moreno
 group; (b)
 $D/\Phi (D)$ is a non-cyclic chief factor of $G$  and if  $D$ is non-abelian,
 then the center, commutator subgroup,
 and Frattini subgroup of $D$  coincide with one another;
 (c)  if  $p>2$, then $D$ is of exponent  $p$,
 for $p =  2$
 the exponent of $D$ is at most 4. From Assertion (b) it follows that
$Q^{G}=G$.

First we show that  $|\Phi (D)|\leq p$. Indeed, assume that  $|\Phi
(D)| >  p$, and let $M$ be a maximal subgroup of $G$ with  $G=DM$  and $Q\leq M$. Then
$M$ is supersoluble, so
$G$ has a 3-maximal subgroup $T$ such that  $Q\leq T$.  Then $T^{G}=G$.
If $T$ is
$S$-quasinormal in $G$, then $G/T_{G}$ is nilpotent by Lemma 2.3(3). Hence
$QT_{G}/T_{G}$ is normal in $G/T_{G}$, which implies that  $QT_{G}=G\leq
M$.   This contradiction shows that $T$ is modular in $G$. Therefore
$G/T_{G}$ is a $P$-group by Lemma 2.1(ii). But then from the $G$-isomorphism
 $$DT_{G}/T_{G}\Phi (D) \simeq D/D\cap T_{G}\Phi (D)=D/\Phi (D)(D\cap
T_{G})=D/\Phi (D)$$ we get that $D/\Phi (D)$ is cyclic. This contradiction
shows that $|\Phi (D)|\leq p$.

Now we show that $|Q|=q$. Assume that $|Q| > q$. Let $M$ be a maximal subgroup of $G$ with
$|G:M|=q$. Then $M$ is supersoluble, so $G$ has a 3-maximal subgroup $T$ such that
 $|G:T|=pq^{2}$.  Then $D\leq T^{G}$  and also  we have $T_{G}\cap D\leq   \Phi (D)$ and $T_{G}
\leq \Phi (D)Q$.     Moreover,  $T$ is not $S$-quasinormal in $G$ since $Q$
 is a Sylow $q$-subgroup
 of $G$ and $|T\cap \Phi (D)| > p$.  Hence $G/T_{G}=(T^{G}/T_{G})\times
(K/T_{G})$ where $T^{G}/T_{G}$ is a non-abelian $P$-group
 of order prime to  $|K/T_{G}|$ by Lemma 2.1(ii), But, clearly, $q$
divides  $|(G/T_{G}):(T^{G}/T_{G})|$, so $Q\leq K$ and hence
$$QT_{G}/T_{G}\leq C_{ G/T_{G}}(T^{G}/T_{G})\leq
C_{G/T_{G}}(DT_{G}/T_{G}\Phi (D)), $$  where  $DT_{G}/T_{G}\Phi (D)$
is $G$-isomorphic to $D/\Phi (D)$. Therefore $|D/\Phi (D)|=p$, a
contradiction.  Hence $|Q|=q$.

Note that if  $|D/\Phi (D)| > p^{2}$ and  $T$ is a subgroup of $D$
with $|D:T|=p^{2}$, then  $T$ is a 3-maximal subgroup of $G$. But
$T$ is neither  $S$-quasinormal nor modular in $G$, a contradiction.
 Hence $|D/\Phi (D)|= p^{2}$.

Finally,  assume that $\Phi (D)\ne 1$,  so $|D|=p^{3}$.
  Assume also  that $p\ne 2$.
 First note that
 since $|Q|=q$, $D$  is a maximal subgroup of $G$. On the other hand, from Assertions (b) and
 (c) we get that some subgroup $T$ of $D$ of order $p$ is not contained in
$\Phi (D)$. It is clear that $T$ is a 3-maximal subgroup of $G$ and  it is
  not $S$-quasinormal in $G$. Hence $T$ is modular in
$G$, so $T^{G}=D$ is a non-abelian $P$-group by Lemma 2.1(ii). This
contradiction shows that $p=2$. Hence $q=3$, since $G/C_{G}(D/\Phi
(D))\simeq Q$ and  $|D/\Phi (D)|=4$.

The proposition is proved.

 {\bf Lemma 3.3.}  {\sl  Suppose that $G$ is soluble and every $ n$-maximal subgroup
 of $ G$
 is either  modular or $S$-quasinormal  in $G$. If $n \leq |\pi (G)|+1$, then
  $G$ is  $\phi$-dispersive for some
 ordering $\phi$ of $\mathbb{P}$.  }

{\bf Proof.}   Suppose that this lemma  is false and let $G$
 be a counterexample of minimal order. Let $N$ be a minimal normal
subgroup of $G$ and $P$  a Sylow $p$-subgroup of $G$  where $p$ divides
$|N|$.  Then $N \leq  P$.

(1) {\sl  $C_{G}(N)=N$ and $G/N$ is strongly supersoluble. Hence $N < P$.}

  Lemma 2.4 implies that the
hypothesis holds for $G/N$. Hence the choice of $G$ implies that $G/N$ is $\phi$-dispersive
for some ordering $\phi$ of $\mathbb{P}$, so $N < P$.  Therefore the choice of $G$ and
Lemma 3.1 imply that   $N\nleq \Phi (G)$.  Hence for some maximal subgroup
 $M$ of $G$ we have  $G=N\rtimes M$.  Then
 $\pi (M)=\pi (G),$ so $G/N\simeq M$ is strongly supersoluble by Theorem 2.12.
Therefore $N$ is a unique minimal normal subgroup of $G$ by
Lemma 2.5. Hence $C_{G}(N)=N$.

(2) $|\pi (G)| > 2$.

 Indeed, assume that  $\pi (G)=\{p, q\}$, and let $Q$
be a Sylow $q$-subgroup of $G$.
  Since $G/N$ is Ore dispersive by Claim (1) and   $P$ is not normal in
$G$,   $NQ/N$ is a normal Sylow subgroup of $G/N$, so for some
normal subgroup $V$ of $G$ we have $N\leq V$ and $|G:V|=p$. Then
$\pi (V)=\pi (G)$. Hence $V$ is   strongly supersoluble by Theorem
2.12. It follows that for the largest prime  $r\in \pi (V)$ a Sylow
$r$-subgroup $R$ of $V$
 is characteristic  in $V$ and so $R$
 is normal in $G$. Hence $r=p$ is the largest prime in $\pi (G)$. Since $M$ is also
 Ore dispersive, a Sylow $p$-subgroup
$M_{p}$ of $M$ is normal in $M$, so it is normal in $G$ since
$N_{G}(M_{p})\nleq M$. But then $NM_{p}$ is a normal Sylow subgroup of
$G$. This contradiction shows that   $|\pi (G)| > 2$.

Take a prime divisor $q$ of the order of $G$ distinct from $p$. Take
a Hall $q'$-subgroup $E$ of $G$, and let $E\leq  W$ where $W$ is a
maximal subgroup of $G$.
 Then $N \leq  E$ and since $G$ is soluble, Lemmas 2.2(4) and 2.3(1) imply that the hypothesis
holds for  $W$. Consequently, the choice of $G$ implies
 that for some prime $t$ dividing
 $|E|$ a Sylow $t$-subgroup $Q$ of $E$ is normal in $E$. Furthermore, since
$C_{G}(N)=N$ we have  $N \leq Q$.
 Hence, $Q$ is a Sylow $p$-subgroup of $E$. It
 is clear also that $Q$  is a Sylow $p$-subgroup of $G$ and $(|G :
N_{G}(Q)|, r) = 1$ for every prime $r \ne  q$. Since $ |\pi(G)| >
2$, it follows that $Q$  is normal in $G$, so $N=Q=P$. This
contradiction completes the proof of the lemma.

In fact, Theorem B is a special case of the following

 {\bf Theorem 3.4.}  {\sl  Suppose that $G$ is soluble and every $ n$-maximal subgroup
 of $ G$
 is either  modular or  $S$-quasinormal   in $G$.  If $n \leq |\pi (G)|+1$ ,
 then $G^{\frak{U}_{s}}$ is a nilpotent Hall subgroup of $G$.  }

 {\bf Proof.}  Suppose that this theorem is false and let $G$
 be a counterexample of minimal order.
Then $G$ is not strongly
supersoluble, so $D=G^{\frak{U}_{s}}\ne 1$.
  By Lemma 3.3, $G$ has a normal  Sylow $p$-subgroup $P$ for some prime $p$
 dividing $|G|$.

(1) {\sl The conclusion of the theorem holds for every quotient
$G/R\ne G/1$} (This directly follows from Lemma 2.4).

(2) {\sl $D$ is nilpotent. }

Assume that this is false. Then,  since $G^{\frak{U}_{s}}\leq G'$,  $G$ is not supersoluble.

 Let $R$ be a minimal normal subgroup of $G$. By   Claim (1) and
 \cite[2.2.8]{Bal-Ez},
$(G/R)^{\frak{U}_{s}}=DR/D \simeq D/D\cap R$ is nilpotent. If $G$ has a
 minimal normal subgroup $N\ne R$,  then $D/D\cap (R\cap N)$ is nilpotent.
Hence $R$ is a unique minimal normal subgroup of $G$ and, by \cite[A,  13.2]{DH},
 $R\nleq \Phi (G)$.  Therefore $R=C_{G}(R)$ by \cite[A, 15.6]{DH},
 and $G=R\rtimes M$
 for some maximal subgroup $M$ of $G$ with $M_{G}=1$. Then $R=P$ is a Sylow $p$-subgroup of
 $G$ by \cite[A, 13.8]{DH}. It is clear that $M$ is not supersoluble, so $|R| > p$
 since otherwise $M\simeq
 G/R=G/C_{G}(R)$ is cyclic.

 Now let $T$ be any
maximal subgroup of $M$. Then $RT$ is a maximal subgroup of $G$ and
 $|\pi (RT)|=|\pi (G)|$ or
 $|\pi (RT)|=|\pi (G)|-1$. Hence,  by
Lemmas 2.2(4) and 2.3(1),   $RM$ satisfies the
 same assumptions as $G$, with $n-1$ replacing $n$. The choice
 of $G$ implies that $(RT)^{\frak{U}_{s}}\leq F(RT)=R$.  Therefore
 $T\simeq T/(T\cap (RT)^{\frak{U}_{s}})\simeq  (RT)^{\frak{U}_{s}}T/(RT)^{\frak{U}_{s}}$ is
 strongly
  supersoluble.  Hence  $M$ is  a
$\frak{U}$-critical group.

By Lemma  2.6(1), $1 < |\pi (M)|\leq 3$.
 First assume that  $|\pi (M)| =2$, then
$n=4$  by Theorem 2.12 since $M$ is not supersoluble. Hence
 every 3-maximal
 subgroup of $M$ is  either modular or $S$-quasinormal in $G$.  Proposition 3.2
 implies that $M$ either is
a non-supersoliuble group of order $qr^{2}$ for some distinct primes $q$ and $r$,
 or $M=Q\rtimes L$,
 where $Q=C_{M}(Q)$ is a quaternion group of order 8 and $|L|=3$.  Then
$R$ is a 3-maximal subgroup of $G$.  Thus every maximal subgroup of $R$ is
either modular or $S$-quasinormal in $G $  and   so  $|R|=p$ by Lemma 2.10(2), a contradiction.
 Thus
  $|\pi (M)| =3$, so $n=5$ and hence  every 4-maximal
 subgroup of $M$ is  either modular or $S$-quasinormal in $G$.
 Let $|M|=q^{a}r^{b}t^{c}$, where
 $p$,  $r$  and $t$ are primes. If $a + b + c > 4$, then some member $T$  of a
composition series of $M$  is a
non-identity  4-maximal subgroup of $M$ since $G$ is soluble, which is
 impossible by Lemma 2.10. Hence  $a +
b + c = 4$ since $M$ is not supersoluble. Therefore, by Lemma 2.6,
$M=Q\rtimes (L\rtimes T)$, where   $|Q|=q^{2}$, $|L|=r$ and $|T|=t$. Then
 $R$ is a 4-maximal subgroup of $G$, so every maximal subgroup of $R$ is
either modular or $S$-quasinormal in $G$, which is impossible by Lemma 2.10(2).
  This contradiction completes the  proof of Claim (2).

(3)  {\sl $D$ is a Hall subgroup of $G$.}

 Suppose that this is false. Then     $G$ is not strongly
supersoluble. Let $P$ be a  Sylow
$p$-subgroup of $D$ such that $1 < P < G_{p}$, where $G_{p}\in
\text{Syl}_{p}(G)$.

(a)  {\sl     $D=P$ is  a minimal normal subgroup of $G$. }

Let $R$ be a minimal normal subgroup of $G$ contained in $D$.
 Then    $R$ is a $q$-group    for some prime
$q$. Moreover, $D/R=(G/R)^{\mathfrak{U}_{s}}$  is a Hall subgroup of
$G/R$ by Claim (1) and \cite[2.2.8]{Bal-Ez}.  Suppose that  $PR/R
\ne 1$. Then  $PR/R \in \text{Syl}_{p}(G/R)$. If $q\ne p$, then
$P \in \text{Syl}_{p}(G)$. This contradicts the fact that $P <
G_{p}$. Hence $q=p$, so $R\leq P$ and therefore $P/R \in
\text{Syl}_{p}(G/R)$ and we again get  that $P \in
\text{Syl}_{p}(G)$. This contradiction shows that  $PR/R=1$, which
implies that
  $R=P$  is the unique minimal normal subgroup of $G$ contained in $D$. Since $D$ is
 nilpotent by Claim  (2),
 a $p'$-complement $E$ of $D$ is characteristic in
$D$ and so it is normal in $G$. Hence $E=1$, which implies that $R=D=P$.

(b) {\sl $D\nleq \Phi (G)$.    Hence for some maximal subgroup
 $M$ of $G$ we have $G=D\rtimes M$  }  (This follows  from  Lemma 2.5 since $G$
 is not strongly supersoluble).

(c) {\sl If $G$ has a minimal normal subgroup $L\ne D$, then $G_{p}=D\times L$.
  Hence $O_{p'}(G)=1$.}

Indeed, $DL/L\simeq D$ is a Hall subgroup of $G/L$ by Claim (1).
Hence  $G_{p}L/L=RL/L$, so $G_{p}=D\times (L\cap G_{p})$. But $D <
G_{p}$, so $(L\cap G_{p})$  is a non-trivial subgroup of $L$. Since
$G$ is soluble, it follows that $L$ is a $p$-group and so
$G_{p}=D\times L$.
 Thus  $O_{p'}(G)=1$.

(d)  $\Phi (G_{p})=  1$.

Suppose that
 $\Phi=\Phi (G_{p})\ne  1$. Then, since $G_{p}$ is normal in $G$ by
Claim (c),  $\Phi$  is normal in $G$ and so we can take a minimal normal subgroup $L$
 of $G$ contained in $\Phi$. But then  $G_{p}=D\times L=D$ by Claim
(c), a contradiction.  Hence we have (d).

{\sl Final contradiction for (3).} Claim (d) implies that $G_{p}$ is
an
 elementary abelian normal subgroup of $G$.
By Maschke's
theorem $G_{p}=N_1\times N_2$ is the direct product of some
 minimal normal subgroups   of $G$.  Claim (a) implies that $N_{1} < G_{p}$.
Let $M=N_{2}E$, where $E$ is a complement to $G_{p}$ in $G$. Then
$M$ is a maximal subgroup of $G$ and $\pi (M)=\pi (G)$.  On the
other hand, every $(n-1)$-maximal subgroup of $M$ is either modular
or $S$-quasinormal in $M$ by hypothesis and Lemmas 2.2(4) and
2.3(1). Thus $M\simeq G/N_{1}$ is strongly supersoluble by Theorem
2.12. Similarly we get that  $G/N_{2}$  is strongly supersoluble.
Hence $G\simeq G/N_{1}\cap N_{2}$ is strongly supersoluble  by Lemma
2.5. This contradiction
 shows that $D=G^{\mathfrak {U}_{s}}$
  a Hall subgroup of $G$.

The proof of the theorem is complete.

\section{Final remarks}

1.  Some preliminary results are of independent significance because they
generalize some known results.

From  Proposition 2.11 we get the following

{\bf Corollary 4.1} (Schmidt \cite{1-3}).   {\sl  If Every
  2-maximal subgroup $M$ of $G$  is modular, then   $G$ is nearly nilpotent. }

{\bf Corollary 4.2.}   {\sl  If
    every    2-maximal subgroup of $G$
 is  $S$-quasinormal in $G$, then $G$ is  nearly nilpotent. }

{\bf Corollary 4.3}   (Agrawal \cite{Agr}). {\sl  If
    every    2-maximal subgroup of $G$
 is  $S$-quasinormal in $G$, then $G$ is  supersoluble. }

From Proposition 3.2 we get the following known result.

{\bf Corollary 4.4} (Schmidt \cite{1-3}).  {\sl  If every
  3-maximal subgroup $M$ of $G$  is modular in $G$ and $G$ is not supersoluble, then  either
    $G$ is
a group of order $pq^{2}$ for some distinct primes $p$ and $q$ or $G=Q\rtimes P$,
 where $Q=C_{G}(Q)$ is a quaternion group of order 8 and $|P|=3$. }

2. In closing note that the restrictions on $|\pi (G)|$  in Theorems A  and
 B cannot be weakened. Indeed, for
Theorem A this follows from  the example of the alternating group $A_{4}$ of degree 4. For
Theorem B this follows from  the example of the  $A_{4}\times C_{2}$, where
$C_{2}$ is a group of order 2.

\end{document}